\documentclass[letterpaper, 9 pt, conference]{ieeeconf}
\IEEEoverridecommandlockouts    
\usepackage[utf8]{inputenc}
\usepackage[T1]{fontenc}
\usepackage{cite}
\usepackage{soul}
\usepackage{amsmath,amssymb,amsfonts}

\usepackage{amsthm}
\usepackage{algorithmic}
\usepackage{algorithm}
\usepackage{graphicx}
\usepackage{graphics}
\usepackage{subfigure}
\usepackage{caption}
\usepackage{subcaption}
\usepackage{booktabs}
\usepackage{multirow}
\usepackage{tcolorbox}
\usepackage{url}
\usepackage{hyperref}
\usepackage[usenames, dvipsnames]{xcolor}
\usepackage{color}
\usepackage{textcomp}
\usepackage{mathtools, cuted}
\usepackage{breqn}
\usepackage{dsfont}
\usepackage{soul}
\usepackage{comment}
\usepackage{lipsum}
\usepackage{stfloats}
\usepackage{float}

\newtheorem{proposition}{Proposition}

\newtheorem{property}{Property}

\newtheorem{definition}{Definition}[section]

\newcommand{\mC}{{\mathbb C}}

\newcommand{\mR}{{\mathbb R}}

\newcommand{\bF}{{\mathbf F}}
\newcommand{\bk}{{\mathbf k}}

\newcommand{\bA}{{\mathbf A}}

\newcommand{\bw}{{\mathbf w}}

\newcommand{\bZ}{{\mathbb Z}}
\newcommand{\cS}{{\mathcal S}}

\newcommand{\cP}{{\mathcal P}}
\newcommand{\cX}{{\mathcal X}}
\newcommand{\cF}{{\mathcal F}}
\newcommand{\cK}{{\mathcal K}}
\newcommand{\cC}{{\mathcal C}}
\newcommand{\cM}{{\mathcal M}}
\newcommand{\mU}{{\mathbb U}}

\newcommand{\bx}{{\mathbf x}}

\newcommand{\bs}{{\mathbf s}}
\newcommand{\bff}{{\mathbf f}}

\title{Spectral Koopman Approach for Reconstructing  State-space Geometry of Cislunar Restricted 3-Body Problem}

\author{Subhrajit Sinha, Sriram S. K. S. Narayanan, Raktim Bhattacharya, and Umesh Vaidya\\
\thanks{S. Sinha is with Pacific Northwest National Laboratory, Richland, WA, USA - 99354. S. Narayanan and U. Vaidya are with Clemson University, Clemson, SC, USA - 29634. R. Bhattacharya is with Texas A\& M University, College Station, TX, USA - 88840. \small{\tt email: subhrajit.sinha@pnnl.gov@pnnl.gov}}
}

\date{}

\begin{document}

\maketitle

\begin{abstract}
In this work, we propose a novel approach, based on the path integral formulation of Koopman spectrum, to discover the phase-space geometry of the planar Cislunar Restricted 3 Body Problem (CR3BP). In contrast to existing techniques, which use trajectory-based (usually) local analysis, we leverage the Koopman operator framework, which generates a global linear \emph{representation} of the system, to reconstruct the global phase space geometry of the CR3BP. In particular, we compute the principal eigenfunctions of the Koopman operator via the path integral approach and show how the zero level curves of these eigenfunctions encode the phase space characteristics of the planar CR3BP.

\end{abstract}

\section{Introduction}

The cislunar environment, defined as the vast region encompassing Earth, the Moon, and the space between them, has emerged as a domain of renewed strategic and scientific interest \cite{holzinger2021primer, whitley2016options}. As national space agencies and commercial entities increasingly target lunar and cislunar operations, understanding the complex dynamical systems governing this region has become essential for mission design and execution. Unlike near-Earth space, cislunar dynamics are shaped by the gravitational interplay of Earth, Moon, and Sun, giving rise to nonlinear phenomena such as libration points, heteroclinic connections, invariant manifolds, and chaotic transport.

From a mission design perspective, spacecraft operations in cislunar space present challenges distinct from traditional orbital regimes. The multi-body gravitational environment leads to highly nonlinear dynamics and sensitive dependence on initial conditions, causing rapid degradation in trajectory prediction accuracy and requiring frequent navigation updates and station-keeping maneuvers \cite{Gomez2001}. These challenges are compounded by sparse tracking infrastructure, high-fidelity perturbation models, and the non-convexity of transfer design, all of which make real-time trajectory optimization particularly difficult.

The current state of the art in cislunar trajectory analysis relies primarily on dynamical systems theory and numerical optimization. The Circular Restricted Three-Body Problem (CR3BP), formulated as a Hamiltonian system with a conserved Jacobi integral~\cite{Szebehely1967,Koon2011}, provides the foundational model, while higher-fidelity extensions such as the Bicircular Four-Body Problem (BCR4BP) and Elliptic Restricted Three-Body Problem (ER3BP) incorporate solar perturbations and orbital eccentricity. Numerical tools including differential correction, pseudo-arclength continuation, multiple shooting, and invariant manifold computation remain central to trajectory design \cite{Howell1984,Gomez2001,Betts1998,Richardson1980,Jorba1999}. Despite their maturity, these approaches face three persistent limitations: fidelity mismatches across models, costly recomputation under navigation updates, and poor scalability of Monte Carlo methods for long-horizon uncertainty propagation.

An alternative approach is to analyze cislunar dynamics through the lens of operator theory, which shifts attention from individual trajectories to the evolution of observables and densities \cite{Lasota}. In this framework, the Koopman operator is especially attractive because it provides a \emph{global} linear description of an otherwise nonlinear dynamical system, allowing nonlinear behavior to be characterized through the spectral properties of a linear operator \cite{Lasota,budivsic2012applied,mauroy2020koopman}. This is particularly appealing in cislunar dynamics, where phase space is organized by strongly nonlinear structures such as invariant manifolds, resonances, and transport channels that are poorly captured by local linearization.

Most existing Koopman-based methods approximate the operator in a finite-dimensional subspace using techniques such as Extended Dynamic Mode Decomposition (EDMD) \cite{williams2015data}. While effective in some settings, these approaches often require careful basis design, substantial data, and can struggle to accurately recover the underlying spectrum \cite{sinha2020robust}. More recently, operator-theoretic ideas have begun to appear in astrodynamics applications, including approximate orbital solutions and local analysis near three-body libration points \cite{arnas2021approximate,servadio2022dynamics}. However, many of these efforts remain either local in scope or dependent on approximation architectures that obscure the intrinsic global geometry of the phase space.

In contrast, the path integral framework introduced in~\cite{deka2023path} offers a direct route for computing Koopman spectral objects without relying on ad hoc dictionary selection. Since the principal eigenfunctions encode invariant and transport structures of the underlying dynamics, they provide a natural mechanism for uncovering the global phase-space geometry relevant to cislunar motion, while also supporting downstream tasks such as uncertainty propagation, reachability, and control \cite{vaidya2025koopman,vaidya2022spectral,deka2024extensions,umathe2022reach} and in this work we demonstrate how the phase-space geometry of a CR3BP can be recovered from time-series data using the path integral formulation.

\section{Preliminaries}

\noindent In this preliminary section, we provide a brief overview of existing results on the spectral theory of the Koopman operator and the path-integral approach for the computation of Koopman eigenfunctions. For more details on this topic, refer to \cite{mezic2020spectrum,mezic2020koopman}.

\subsection{Spectral Theory of Koopman Operator}
Consider the dynamical system
\begin{align}
    \dot \bx={\bff}(\bx),\label{odesys}
\end{align}
defined on a state space ${\cal X}\subseteq \mR^p$.  The vector field $\bff$ is assumed to be a smooth function. Let $\cF\subseteq \cC^0$ be the function space of observable $\psi: \cX\to \mC$.
We have the following definitions for the Koopman operator and its spectrum \cite{mezic2020spectrum}. 
\begin{definition}[Koopman Operator] The family of Koopman  operators $\mathbb{U}_t:\cF\to \cF$ corresponding to ~\eqref{odesys} is defined as 
\begin{eqnarray}[\mathbb{U}_t \psi](\bx)=\psi(s_t(\bx)). \label{koopman_operator}
\end{eqnarray}
If in addition $\psi$ is continuously differentiable, then $f(\bx,t):=[\mU_t \psi ](\bx)$ satisfies a partial differential equation \cite{Lasota} 
\begin{align}
\frac{\partial f}{\partial t}=\frac{\partial f}{\partial \bx} \bff=: \cK_\bff f \label{Koopmanpde}
\end{align}
with the initial condition $f(\bx,0)=\psi(\bx)$. The operator $\cK_\bff$ is the infinitesimal generator of $\mU_t$ i.e.,
\begin{eqnarray}
{\cal K}_{\bff} \psi=\lim_{t\to 0}\frac{(\mathbb{U}_t-I)\psi}{t}. \label{K_generator}
\end{eqnarray}
\end{definition}
It is easy to check that each $\mU_t$ is a linear operator on the space of functions, $\cF$.  
\begin{definition}\label{definition_koopmanspectrum}[Eigenvalues and Eigenfunctions of Koopman] A function $\psi_\lambda(\bx)$, assumed to be at least $\cC^1$,  is said to be an eigenfunction of the Koopman operator associated with eigenvalue $\lambda$ if
\begin{eqnarray}
[\mU_t \psi_\lambda](\bx)=e^{\lambda t}\psi_\lambda(\bx)\label{eig_koopman}.
\end{eqnarray}
Using the Koopman generator, the (\ref{eig_koopman}) can be written as 
\begin{align}
    \frac{\partial \psi_\lambda}{\partial \bx}{\bff}=\lambda \psi_\lambda\label{eig_koopmang}.
\end{align}
\end{definition}
With the hyperbolicity assumption on the equilibrium point of the system (\ref{odesys}), the principal eigenfunctions are well-defined \cite{mezic2020spectrum}.

The principal eigenfunctions can be used as a change of coordinates in the linear representation of a nonlinear system and draw a connection to the famous Hartman-Grobman theorem on linearization and Poincaré normal form \cite{lan2013linearization,mezic2020spectrum,arnold2012geometrical}. 
The principal eigenfunctions will be defined over a proper subset $\cP$ of the state space $\cX$ (called subdomain eigenfunctions) or over the entire $\cX$ \cite[Lemma 5.1, Corollary 5.1, 5.2, and 5.8]{mezic2020spectrum}. For a system with a stable equilibrium point, the principal eigenfunctions are well-defined in the entire domain of attraction of the equilibrium point.

\begin{property}\label{property_algebra} [Algebra structure of eigenfunction under product] Let $\varphi_{\lambda_1}, \varphi_{\lambda_2}\in \cF$ be the eigenfunctions with eigenvalues $\lambda_1$ and $\lambda_2$ then $\varphi_{\lambda_1}^{k_1}\varphi_{\lambda_2}^{k_2}$ is also an eigenfunction for $k_1,k_2\in {\mathbb Z}^+$ with eigenvalues $k_1\lambda_1+k_2\lambda_2$.
\end{property}

\begin{definition}[Koopman Mode Decomposition (KMD)]\label{definition_kmd} Consider a time dependent function $\psi(\bx)\in \cF$ , we can decompose the function using Koopman eigenfunctions as 
\begin{align}
\psi(\bx)=\sum_{\bk\in \bZ^n} c_{\bk}\varphi_{\lambda_1}^{k_1}(\bx)\ldots \varphi_{\lambda_n}^{k_n}(\bx)
\end{align}
The coefficients $c_\bk$ are the Koopman modes and correspond to the projection of the function $\psi$ along the eigenfunctions $\varphi_{\lambda_1}^{k_1}(\bx)\ldots \varphi_{\lambda_n}^{k_n}(\bx)$. 
\end{definition}

\subsection{Computing the Principal Spectrum using Path Integrals}
{\color{black} 
We begin by decomposing the system \eqref{odesys} as the sum of linear and nonlinear components:
\begin{align}
\dot \bx=\bff(\bx)=\bA\bx+\bF_n(\bx)\nonumber
\end{align}
where $\bA=\frac{\partial \bff}{\partial \bx}(0)$ denotes the Jacobian of the vector field at the origin, and $\bF_n(\bx)=\bff(\bx)-\bA\bx$ represents the purely nonlinear part.

Similarly, the principal eigenfunction of the Koopman operator associated with eigenvalue $\lambda$ (and the corresponding left eigenvector $\bw_\lambda$) can also be separated into linear and nonlinear terms~\cite{mezic2020spectrum}:
\begin{align}
\varphi_\lambda(\bx)=\bw_\lambda^\top \bx+h_\lambda (\bx)\label{eig_decomposition}
\end{align}

By substituting~(\ref{eig_decomposition}) into~(\ref{eig_koopmang}) and separating linear and nonlinear terms, we find that the nonlinear part $h_\lambda$ satisfies the linear PDE
\begin{align}
\frac{\partial h_\lambda}{\partial \bx}\bff(\bx)-\lambda h_\lambda(\bx)+\bw_\lambda^\top \bF_n(\bx)=0. \label{linearpde}
\end{align}

The path-integral solution of $h_\lambda$ is obtained by applying the method of characteristics to the PDE above. Specifically, the solution of~(\ref{linearpde}) can be expressed as~\cite{deka2023path}
\begin{align}
h_\lambda(\bx)=e^{-\lambda t}h_\lambda(\bs_t(\bx))+\int_0^t e^{-\lambda \tau}\bw_\lambda^\top \bF_n(\bs_\tau(\bx)),d\tau \nonumber
\end{align}
Moreover, if $\displaystyle \lim_{t\to \infty}e^{-\lambda t}\bw_\lambda^\top \bF_n(\bs_t(\bx))=0$, then the explicit expression for $h_\lambda(\bx)$ is obtained as
\begin{align}
h_\lambda(\bx)=\int_0^\infty e^{-\lambda \tau}\bw_\lambda^\top \bF_n(\bs_\tau(\bx))\,d\tau \label{eqn:pde_sol}
\end{align}

In practice, we can approximate the principal eigenfunctions by evaluating the integral in \eqref{eqn:pde_sol} up to a finite time $T$ as follows:
\begin{align}
\hat \varphi_{\lambda_u}(\bx) &= \bw^\top_{\lambda_u}\bx + \int_0^T e^{-\lambda_u t},\bw_{\lambda_u}^\top \bF_n(\bs_t(\bx))\,dt, \label{eqn:eigfn_unstable}\\
\hat \varphi_{\lambda_s}(\bx) &= \bw^\top_{\lambda_s}\bx + \int_0^{T} e^{\lambda_s t},\bw_{\lambda_s}^\top \bF_n(\bs_{-t}(\bx))\,dt. \label{eqn:eigfn_stable}
\end{align}
Note that, for $\lambda_u>0$, the unstable eigenfunction is obtained using~\eqref{eqn:eigfn_unstable} using the forward trajectory $\bs_t(\bx)$. For $\lambda_s<0$, the stable eigenfunction is computed using~\eqref{eqn:eigfn_stable} with the backward trajectory $\bs_{-t}(\bx)$.

}


\section{Cislunar Restricted 3 Body Problem}
In this section, we briefly describe the phase-space geometry of the CR3BP. For details we refer the reader to \cite{Gomez2001,arnold_celestial_mechanics}.
\subsection{The Model}
We consider two masses $m_1$ and $m_2$ moving on circular orbits about their barycenter. With the mass parameter $\mu = \frac{m_2}{m_1+m_2}$, in rotating frame, the location of the primaries $m_1$ and $m_2$ are given by $(-\mu, 0, 0)$ and $(1-\mu,0,0)$ respectively. Now, for a test particle, with position $(x,y,z)$, the effective potential is given by
\begin{eqnarray}
\Omega(x,y,z) = \tfrac{1}{2}(x^2+y^2) + \frac{1-\mu}{r_1} + \frac{\mu}{r_2},
\end{eqnarray}
where 
\begin{eqnarray*}
r_1 = \sqrt{(x+\mu)^2 + y^2 + z^2}; \; r_2 = \sqrt{(x-1+\mu)^2 + y^2 + z^2},
\end{eqnarray*}
 are the distances of the test particle from the primaries. 

With this, the equations of motion are given by
\begin{equation}\label{CR3BP_6D}
\begin{aligned}
\ddot{x} - 2\dot{y} &= x - \frac{(1-\mu)(x+\mu)}{r_1^3}
                        - \frac{\mu(x-1+\mu)}{r_2^3}, \\
\ddot{y} + 2\dot{x} &= y - \frac{(1-\mu)y}{r_1^3}
                        - \frac{\mu y}{r_2^3},\;\ddot{z} = - \frac{(1-\mu)z}{r_1^3}
             - \frac{\mu z}{r_2^3}.
\end{aligned}
\end{equation}


Although \eqref{CR3BP_6D} captures the full dynamics of the system \cite{Szebehely1967,murray1999}, in many applications, however, the dynamics is confined to the plane of motion of the primaries. By setting $z=\dot{z}=0$, the system reduces to a four-dimensional \emph{planar CR3BP} in $\bx = (x, \; y,\; \dot{x}, \;\dot{y})$, which is algebraically simpler yet retains the essential nonlinear features needed to study Lyapunov orbits and related planar phenomena \cite{Szebehely1967,Gomez2001}, and in this work, we considered this 4D planar model. 

With this, the equations of motion for the planar 4D CR3BP are
\begin{equation}\label{CR3BP_4D}
\begin{aligned}
\dot{x}   &= v_x, \quad \dot{y} = v_y, \\
\dot{v}_x &= 2v_y + x - \frac{(1-\mu)(x+\mu)}{r_1^3} - \frac{\mu(x-1+\mu)}{r_2^3}, \\
\dot{v}_y &= -2v_x + y - \frac{(1-\mu)y}{r_1^3} - \frac{\mu y}{r_2^3},
\end{aligned}
\end{equation}
where $r_1$ and $r_2$ denote the distances from the infinitesimal mass to the two primaries.  
The system admits the conserved Jacobi integral,
\begin{equation}\label{Jacobi_constant}
C = x^2 + y^2 + 2\left( \frac{1-\mu}{r_1} + \frac{\mu}{r_2} \right)
- \left(v_x^2 + v_y^2\right) = -2H,
\end{equation}
where $H$ is the Hamiltonian of the system. For additional details, see~\cite{Szebehely1967,Gomez2001}.

The planar CR3BP exhibits a rich set of nonlinear dynamical behaviors arising from the interplay between the gravitational attraction of the primaries and the Coriolis forces in the rotating frame. The system admits equilibrium points (the Lagrange points), around which families of periodic and quasi-periodic orbits exist. These orbits undergo bifurcations as system parameters, such as the Jacobi constant, vary, giving rise to stable and unstable manifolds that structure transport and chaotic motion in phase space. The coexistence of regular (near-integrable) dynamics and strongly chaotic trajectories highlights the fundamentally nonlinear character of the CR3BP and makes it a canonical model for studying complex celestial mechanics phenomena such as resonances, orbit transfers, and capture.

\subsection{Phase Space Geometry of the Planar CR3BP}

The planar circular restricted three--body problem (CR3BP) provides a canonical model to study transport dynamics in celestial mechanics. A test particle of negligible mass moves under the gravitational influence of two primaries (Earth and Moon) in a rotating reference frame where the primaries remain fixed. The system admits the Jacobi integral, $C(x,y,\dot{x},\dot{y}) = 2U(x,y) - \big(\dot{x}^2 + \dot{y}^2 \big)$,
where $U(x,y)$ is the effective potential including centrifugal and gravitational contributions. Since the Jacobi constant is the only global integral of motion, the dynamics are restricted to the three--dimensional manifold
\[
\mathcal{M}_{C^\star} = \big\{ (x,y,\dot{x},\dot{y}) \in \mathbb{R}^4 : C(x,y,\dot{x},\dot{y}) = C^\star \big\},
\]
which is a codimension--1 hypersurface of the full four--dimensional phase space.

\begin{figure}[htp!]
\centering
\includegraphics[width=0.45\linewidth]{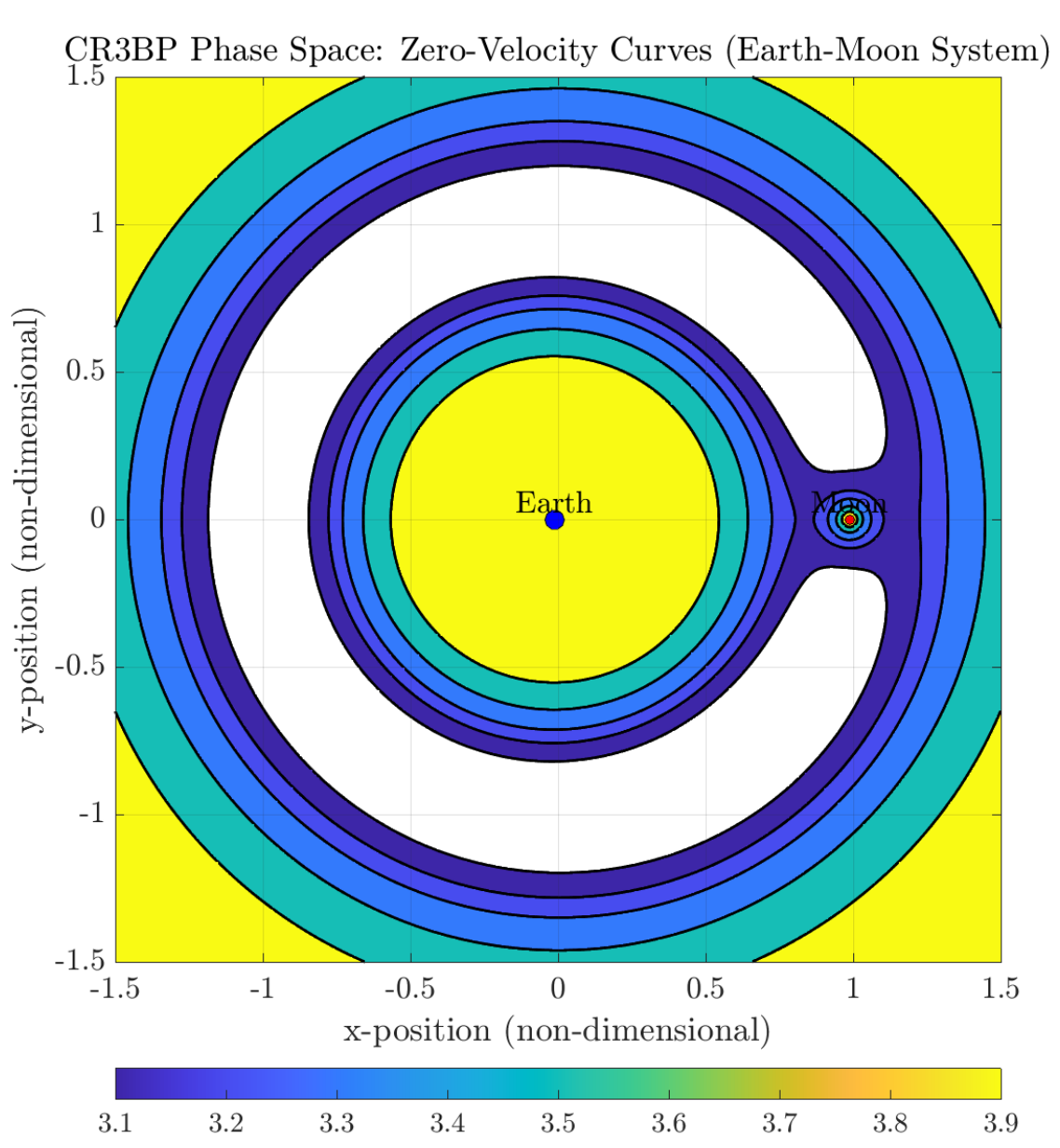}
\caption{Zero--velocity curves for the Earth--Moon CR3BP. The shaded regions indicate forbidden zones, while the corridors represent energetically allowed regions of motion. Bottlenecks near $L_1$ and $L_2$ serve as gateways for transport.}\label{fig_Jacobi_shell}
\end{figure}

\noindent\textbf{Zero--velocity curves in the configuration plane.}  
Fig. \ref{fig_Jacobi_shell} illustrates the classical zero--velocity curves (ZVCs) in the $(x,y)$ configuration plane. These are obtained by setting the velocity terms $\dot{x}=\dot{y}=0$ in the Jacobi integral, yielding the condition $C(x,y,0,0) = 2U(x,y)$. The shaded yellow regions correspond to forbidden zones where $2U(x,y) < C^\star$, since no real velocity values can satisfy the integral. The dark blue corridors represent allowed regions of motion, whose geometry depends strongly on the chosen Jacobi constant. The figure shows how the Hill’s regions separate into disconnected lobes around the Earth and Moon, and how narrow channels open near the equilibrium points $L_1$ and $L_2$. These bottlenecks are gateways for transport: a trajectory with the correct energy may cross between the Earth region, the Moon region, and the exterior domain only through these corridors. Thus, the ZVCs provide a global two--dimensional view of the energetically accessible configuration space.

\begin{figure}[htp!]
\centering
\includegraphics[width=0.55\linewidth]{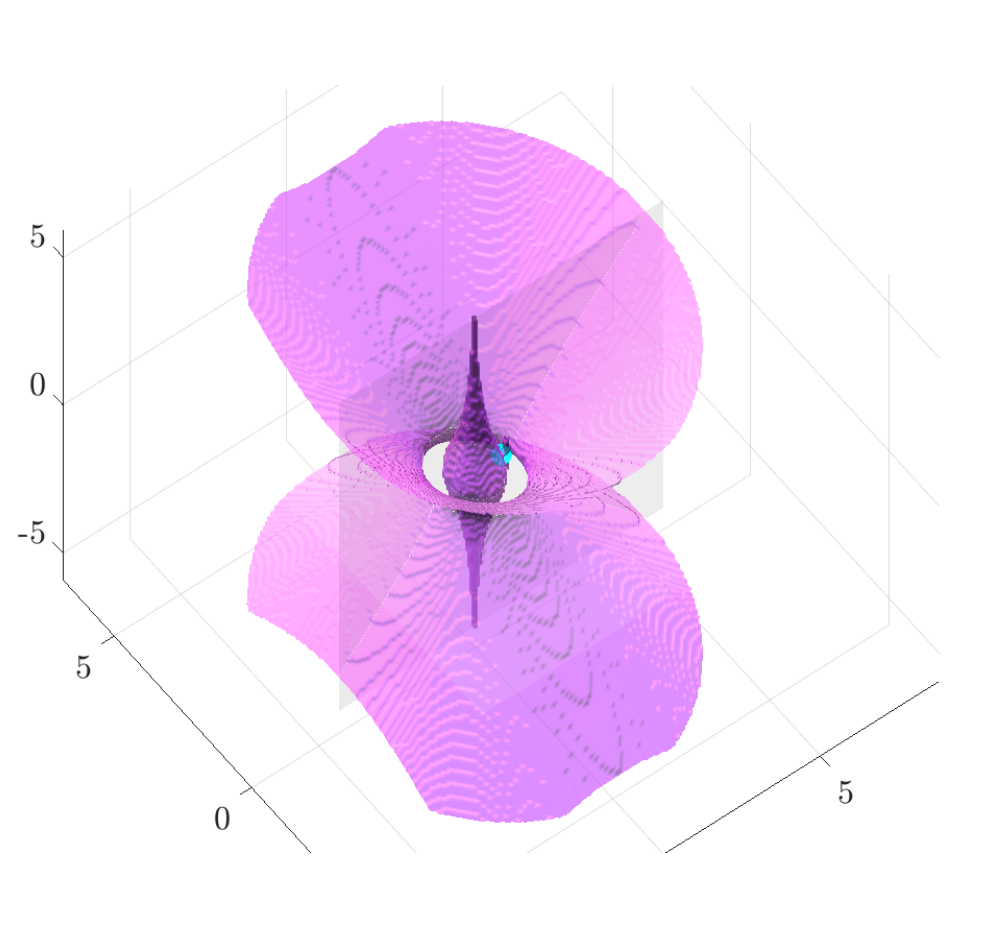}
\caption{Jacobi energy shell near $L_1$ visualized in $(x,y,p_x)$ coordinates.}\label{fig_Jacobi_3D}
\end{figure}

\noindent\textbf{Jacobi energy shell near $L_1$.}  
Fig. \ref{fig_Jacobi_3D} provides a complementary three--dimensional perspective by embedding the Jacobi manifold near the collinear equilibrium $L_1$. We work in coordinates $(x,y,p_x)$, where $p_x = \dot{x}-y$ is the canonical momentum. For each point $(x,y,p_x)$ consistent with $C=C^\star$, the remaining momentum $p_y$ can be recovered from the Jacobi integral, giving two branches $p_y = \pm \sqrt{\,2U(x,y) - C^\star - p_x^2\,}$. These two branches form the upper and lower sheets of the Jacobi energy shell. The plot shows these sheets as semi--transparent surfaces, with the equilibrium point $L_1$ marked in cyan. 

Geometrically, the Jacobi surface is a curved ``energy shell'' in phase space on which all nearby trajectories evolve. Periodic orbits, such as the planar Lyapunov family, appear as closed loops lying entirely on this shell. The invariant manifolds of these periodic orbits, which guide transport between regions, are submanifolds embedded in the Jacobi surface. By rotating the 3D plot, one directly perceives how the local structure near $L_1$ connects the Earth and Moon regions through phase space channels.

\section{Reconstruction of CR3BP Geometry}

The spectral properties of the Koopman operator provide a powerful framework for addressing a wide range of problems relevant to the analysis and design of missions in the cislunar environment. In particular, these properties can be leveraged to identify invariant manifolds, perform uncertainty propagation and reachability analysis \cite{umathe2022reach}, and develop optimal control strategies \cite{vaidya2022spectral,vaidya2025koopman}.

\subsection{Identifying invariant manifolds}

The spectral analysis of the Koopman operator can be used for the identification of invariant manifolds in dynamical systems. In particular, the stable, unstable, and center manifolds can be obtained as joint zero-level curves of the Koopman principal eigenfunctions. We have the following proposition from (Proposition 3.1\cite{mezic2020spectrum}).

\begin{proposition}\label{proposition_mainfolds}
Let $\bx_0$ be the hyperbolic equilibrium point of the system (\ref{odesys}). Let $\lambda_1,\ldots,\lambda_u$ be eigenvalues with positive real part with associated eigenfunctions $\phi_{\lambda_1},\ldots, \phi_{\lambda_u}$ and $\lambda_{u+1},\ldots,\lambda_n$ be eigenvalues with negative real part with associated eigenfunctions $\phi_{\lambda_{u+1}},\ldots, \phi_{\lambda_n}$. Then the joint level set of the eigenfunctions  
\begin{align}
{\cal M}_s=\{\bx\in M: \phi_{\lambda_1}(\bx)=\ldots=\phi_{\lambda_u}(\bx)=0\},    
\end{align}
forms the stable manifold and the joint level set of the eigenfunctions
\begin{align}
{\cal M}_u=\{\bx\in M: \phi_{\lambda_{u+1}}(\bx)=\ldots=\phi_{\lambda_n}(\bx)=0\},
\end{align}
is the unstable manifold of the equilibrium point $\bx_0$. 
\end{proposition}
For center manifolds, we have the following Proposition from (Proposition 5.1\cite{mezic2020spectrum}).

\begin{proposition} Let the system (\ref{odesys}) has $u$ unstable  eigenfunctions $\phi_{\lambda_1},\ldots, \phi_{\lambda_u}$ and $s$ stable  eigenfunctions $\varphi_{\lambda_{u+1}},\ldots, \varphi_{\lambda_{u+s}}$ of the Koopman operator. The joint zero-level set of unstable  eigenfunctions 
\begin{align}
{\cal M}_{cs}=\{\bx\in \cM: \varphi_{\lambda_1}=0,\ldots,\varphi_{\lambda_u}=0 \}
\end{align}
is the center-stable manifold.

\begin{align}
{\cal M}_c=\{\bx\in \cM: \varphi_{\lambda_1}=0,\ldots,\varphi_{\lambda_u}=0, \nonumber \\ 
\qquad \varphi_{\lambda_u+1}=0,\ldots,\varphi_{\lambda_{u+s}}=0, \}
\end{align}

is the center manifold, and 
\begin{align}
{\cal M}_{cu}=\{\bx\in \cM: \varphi_{\lambda_{u+1}}=0,\ldots,\varphi_{\lambda_{u+s}}=0 \}
\end{align}
\end{proposition}

The results from the above proposition can be used for the identification of invariant manifolds of interest in cislunar dynamics. In particular, joint zero-level curves of the Koopman principal eigenfunctions associated with appropriate eigenvalues can be used to identify various invariant manifolds. 

\subsection{Uncertainty propagation and reachability analysis}
The principal eigenfunctions of the Koopman operator can also be used for uncertainty analysis and computation of reachable sets in cislunar dynamics \cite{umathe2022reach}. In particular, let the set of initial conditions with uncertainty be described using the set $\cS_0$ and expressed in terms of a scalar value function, say $g(\bx)\in \cF$, i.e., 
\begin{align}
\cS_0=\{\bx\in {\cal M}: g(\bx)\leq 0\}
\end{align}
The uncertainty in the initial set can be propagated forward or backward in time using the principal eigenfunctions. Let $\cS_t$ be the forward propagation of the set $\cS_0$, i.e., $\cS_t$ be the set of all initial conditions that can be reached from the initial set $\cS_0$. Similarly, $\cS_{-t}$ is the set of all initial conditions that can be steered to the set $\cS_0$ in time $t$. We have following characterization of $\cS_t$ and $\cS_{-t}$ \cite{umathe2022reach}. 
\begin{align}
\cS_t&=\{\bx\in \cM: [\mU_{-t} g](\bx)=:g_f(\bx,t)\leq 0\},\nonumber\\
\cS_{-t}&=\{\bx\in \cM: [\mU_{t} g](\bx)=:g_b(\bx,t)\leq 0\}\label{fb}
\end{align}

Now let $g$ admit the following Koopman mode decomposition in terms of the higher order principal eigenfunctions

\begin{align} g(\bx)=\sum_{\bk\in \mathbb{N}^n}^N\bar g_\bk \prod\limits_{i = 1}^{n} \phi_{\lambda_i}^{k_i}(\bx),
\end{align} $\bk=(k_1,\ldots, k_n)$ and $N$ is the number of eigenfunctions used in the expansion and $\bar g_{\bk}$ are the coefficients used in the expansion. Following (\ref{fb}) and using the definition of eigenfunctions, we can write
\[g_f(\bx,t)=\sum_{\bk\in \mathbb{N}^n}^N\bar g_\bk \prod\limits_{i = 1}^{n} e^{-k_i\lambda_i t} \phi_{\lambda_i}^{k_i}(\bx)\]
\[g_b(\bx,t)=\sum_{\bk\in \mathbb{N}^n}^N\bar g_\bk \prod\limits_{i = 1}^{n} e^{k_i\lambda_i t} \phi_{\lambda_i}^{k_i}(\bx)\]
The above formula for the forward and backward propagation of the initial set is obtained using the algebra property of the principal eigenfunctions (Property \ref{property_algebra}).



\section{Simulation Results}

In the following sections, we demonstrate the effectiveness of the proposed Koopman framework for analyzing planetary systems. Specifically, we employ the spectrum of the Koopman operator, computed via the path-integral approach, to recover fundamental dynamical features of the planar CR3BP.

\subsection{Jacobi Constant and Zero--Velocity Curves}

In the planar circular restricted three--body problem (CR3BP), the motion of a test particle in the rotating frame admits a conserved quantity called the \emph{Jacobi constant} $C$.  
By setting the velocity components $(v_x,v_y)=(0,0)$, one obtains the \emph{zero--velocity curves} (ZVCs), which form the boundaries of energetically allowed motion in the $(x,y)$ plane.  
The regions enclosed by these curves are dynamically accessible; regions outside are forbidden.  

\begin{figure}[htp!]
    \centering
    \includegraphics[width=0.6\linewidth]{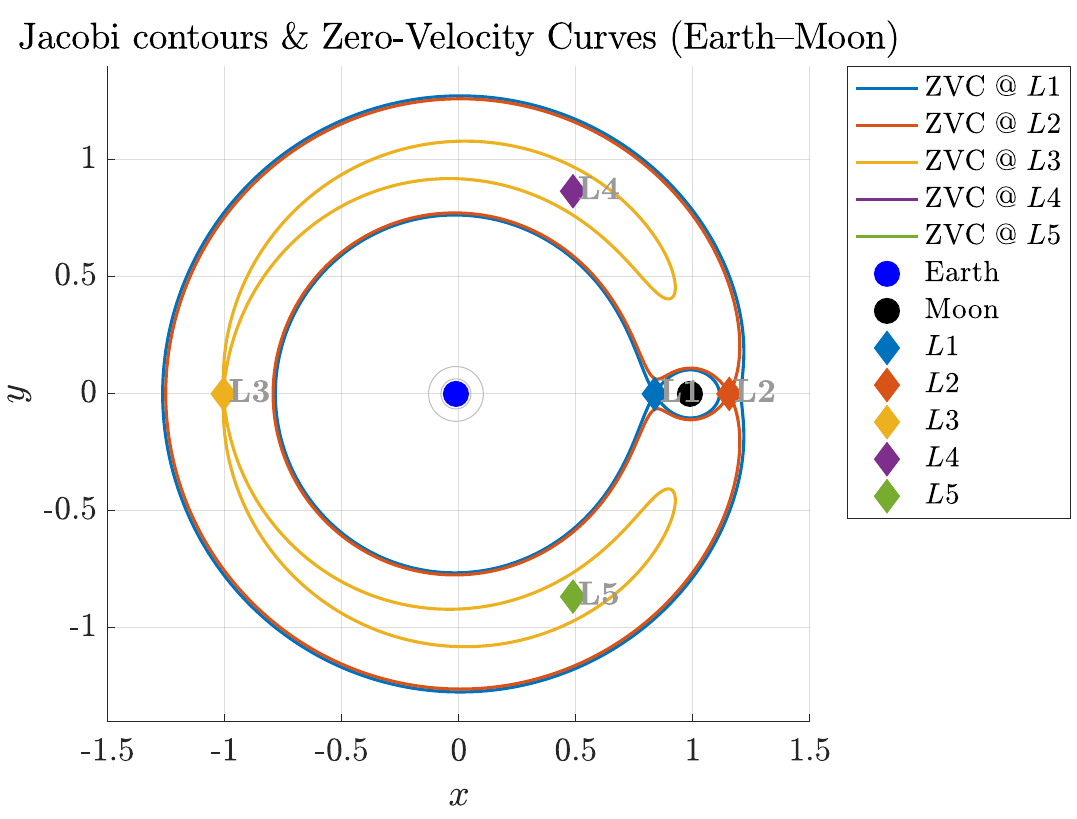}
    \caption{Jacobi contours and zero--velocity curves for the Earth--Moon CR3BP.}\label{fig_jacobi}
\end{figure}

Fig. \ref{fig_jacobi} shows Jacobi contours and highlights the five Lagrange points L$_1$--L$_5$.  
The ``vessel--shaped'' boundaries around the Earth and Moon define \emph{domains of attraction}, within which motion is confined.  
Openings at the collinear points (L$_1$, L$_2$, L$_3$) create \emph{necks} through which transport can occur.  
For example, at certain $C^*$ values, the ZVC near L$_1$ opens, allowing trajectories to leak from Earth’s region toward the Moon.  
These geometric features provide the global scaffold for transport, but the fine structure of actual trajectories is determined by invariant manifolds associated with the equilibria.

\subsection{Koopman Eigenfunctions at L$_1$ and L$_2$}

The principal eigenfunctions of the Koopman operator, $\phi(\bx)$, are computed here using the path--integral formulation given in \eqref{eqn:eigfn_unstable} and \eqref{eqn:eigfn_stable}. Their zero--level sets approximate invariant manifolds:  
unstable ($\phi_u=0$), stable ($\phi_s=0$), and center ($\mathrm{Re}(\phi_c)=0$).  
These functions encode the geometry of expansion, contraction, and oscillation near each equilibrium.

\subsubsection{L$_1$}
\begin{figure}[htp!]
\centering
\subfigure[]{\includegraphics[width=0.4\linewidth]{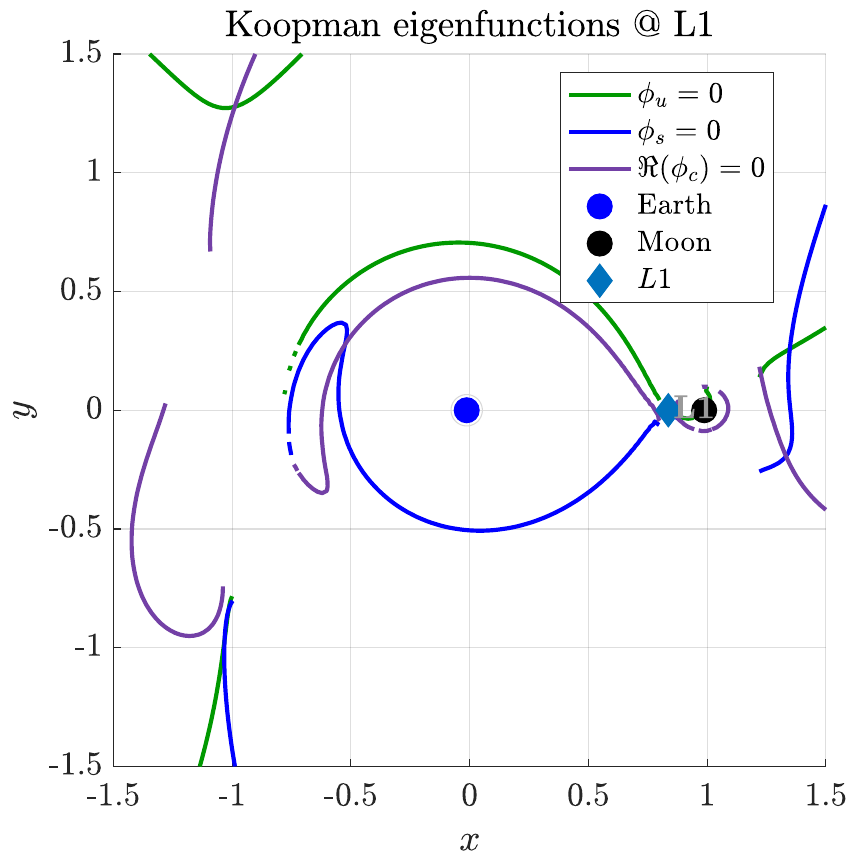}}
\subfigure[]{\includegraphics[width=0.4\linewidth]{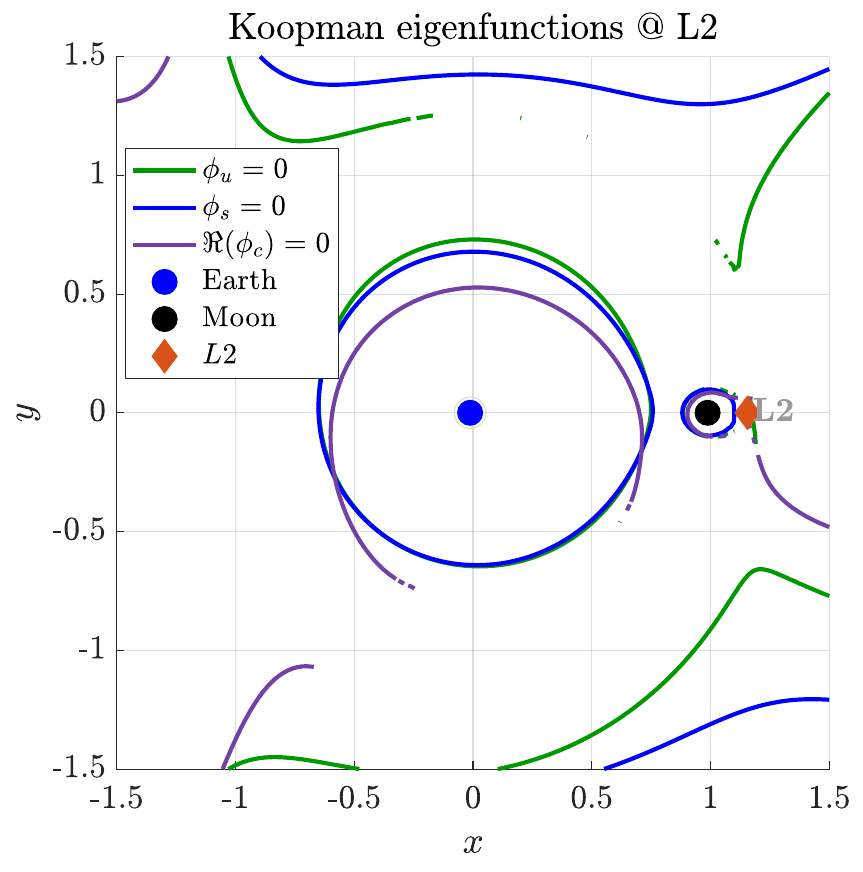}}
\caption{Koopman eigenfunctions at (a) L$_1$ and (b) L$_2$.}\label{fig_L1}
\end{figure}

At L$_1$ (Fig. \ref{fig_L1}(a)), between Earth and Moon, the unstable manifold (green) extends outward through the neck toward the Moon, while the stable manifold (blue) bends inward toward Earth.  
This interplay marks L$_1$ as a transport gateway: trajectories initialized near $\phi_u=0$ can escape Earth’s potential well, while those aligned with $\phi_s=0$ are captured back toward Earth.  
The center eigenfunction (purple) delineates local oscillations, producing loops that encircle the equilibrium.  
Thus the Koopman representation captures both the escape/capture channel and the small oscillatory Lyapunov orbits around L$_1$.  

\subsubsection{L$_2$}

At L$_2$ (Fig. \ref{fig_L1}(b)), just beyond the Moon, the geometry complements that of L$_1$.  
The unstable manifold projects outward into the exterior region, guiding escape trajectories away from the Earth--Moon system.  
Conversely, the stable manifold funnels trajectories inward from the exterior, representing possible capture pathways.  
The center eigenfunction again identifies oscillatory modes around equilibrium, corresponding to small Lyapunov orbits at L$_2$.  
Together, L$_1$ and L$_2$ form a coupled transport system: the unstable manifold of one can align with the stable manifold of the other, creating a channel across the Moon’s neck.

\subsection{Koopman Eigenfunctions at L$_3$, L$_4$, and L$_5$}

\subsubsection{L$_3$}
\begin{figure}[htp!]
\centering
\subfigure[]{\includegraphics[width=0.4\linewidth]{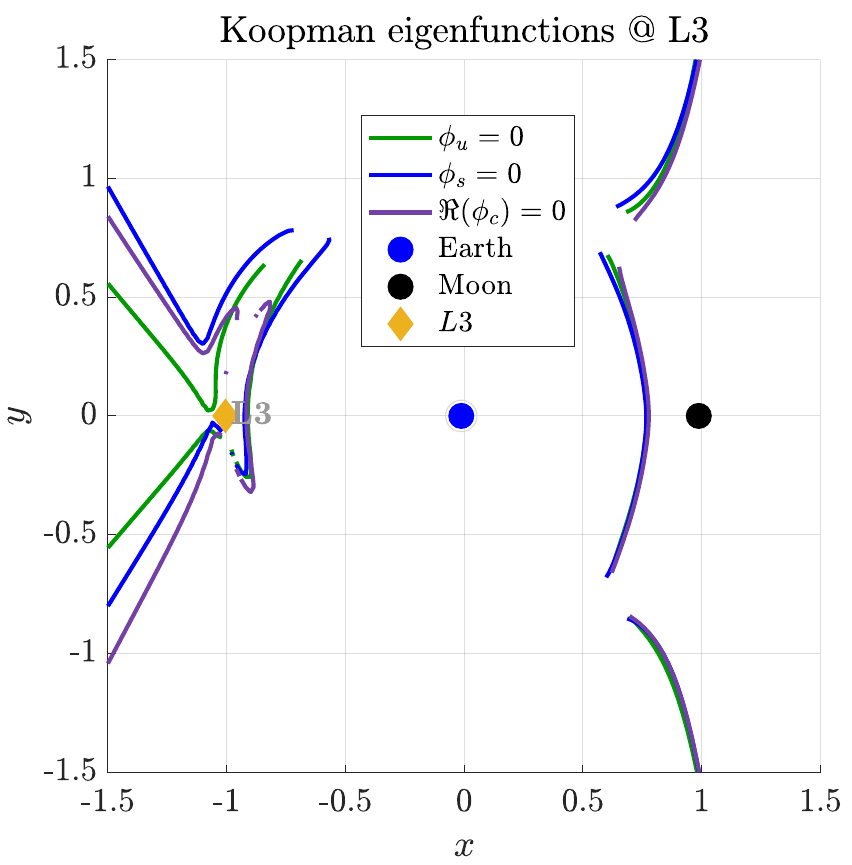}}
\subfigure[]{\includegraphics[width=0.4\linewidth]{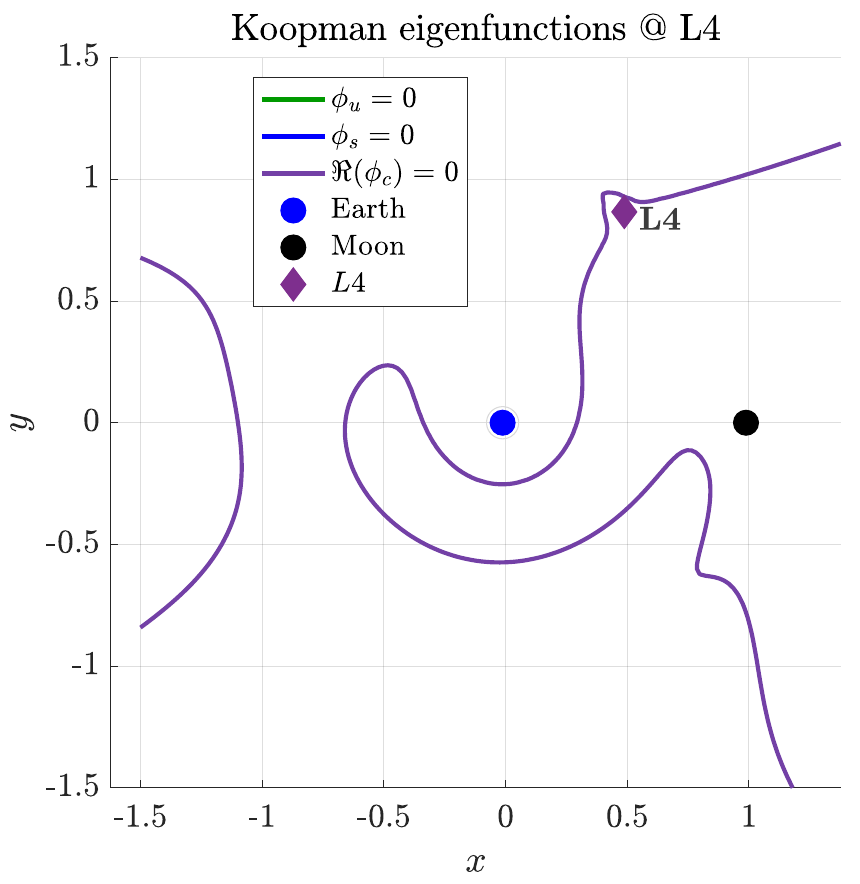}}
\subfigure[]{\includegraphics[width=0.4\linewidth]{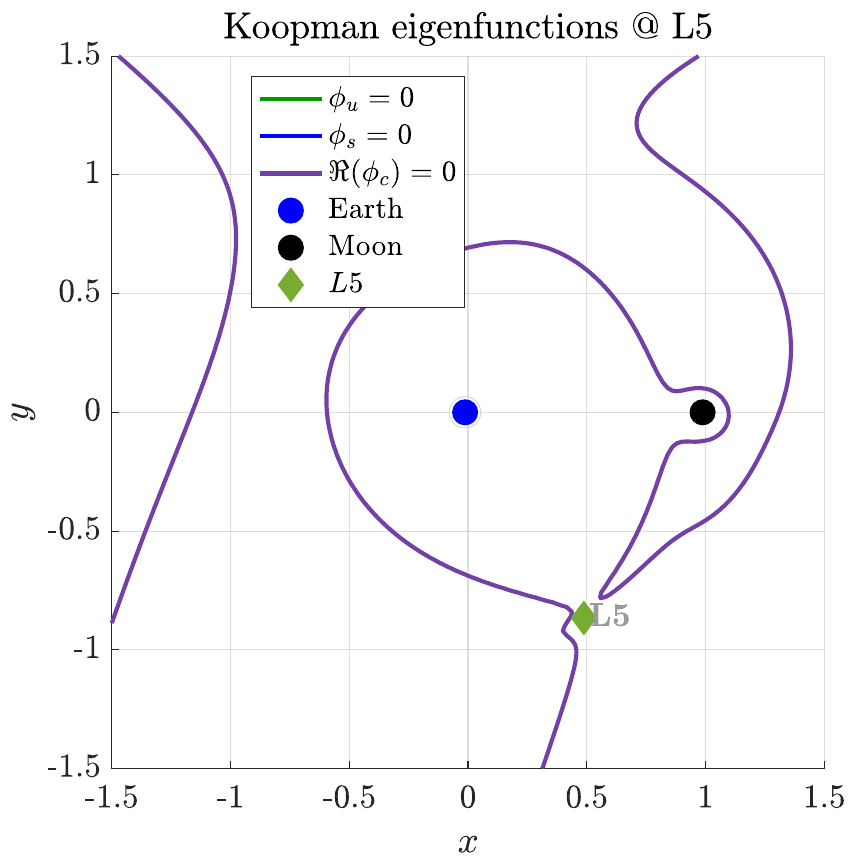}}
\caption{Koopman eigenfunctions at (a) L$_3$. (b) L$_4$. (c) L$_5$.}\label{fig_L3}
\end{figure}

L$_3$ (Fig. \ref{fig_L3}(a)), on the far side of Earth opposite the Moon, exhibits nearly symmetric manifolds about the $y$--axis.  
The unstable and stable eigenfunctions open outward on opposite sides, indicating that motion near L$_3$ is fragile and tends to leak toward the external region rather than funnel into a bounded corridor.  
The center eigenfunction traces oscillatory loops surrounding Earth.  
This highlights that while L$_3$ is mathematically an equilibrium, its manifolds do not support strong transport channels comparable to L$_1$ and L$_2$.  

\subsubsection{L$_4$ and L$_5$}


At the triangular Lagrange points L$_4$ and L$_5$ (Figs. \ref{fig_L3}(b) and \ref{fig_L3})(c), leading and trailing the Moon by $60^\circ$, the eigenfunctions show a very different character.  
Here the center eigenfunction dominates, producing closed loops that encircle the equilibrium.  
These loops correspond to tadpole and Trojan orbits, oscillatory regions that are stable in the planar CR3BP.  
The stable and unstable manifolds are weak, reflecting the fact that these points are linearly stable.  
Thus, the Koopman eigenfunctions successfully reproduce the qualitative difference between collinear (unstable) and triangular (stable) equilibria.

\subsection{Global Picture}

\begin{figure}[htp!]
\centering
\subfigure[]{\includegraphics[width=0.4\linewidth]{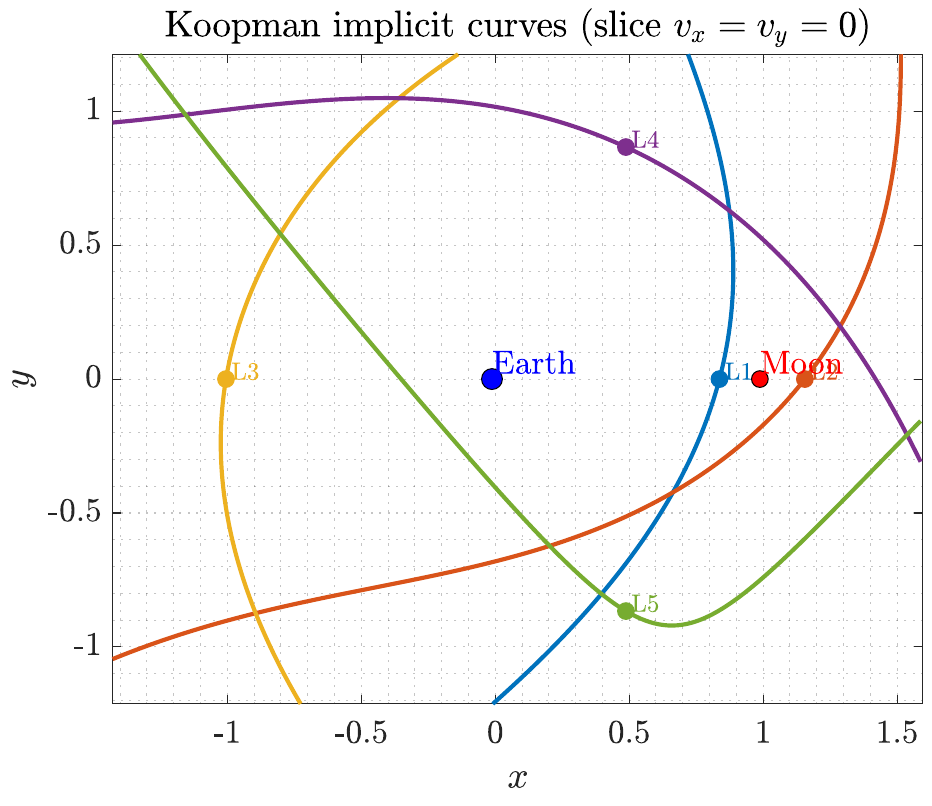}}
\subfigure[]{\includegraphics[width=0.47\linewidth]{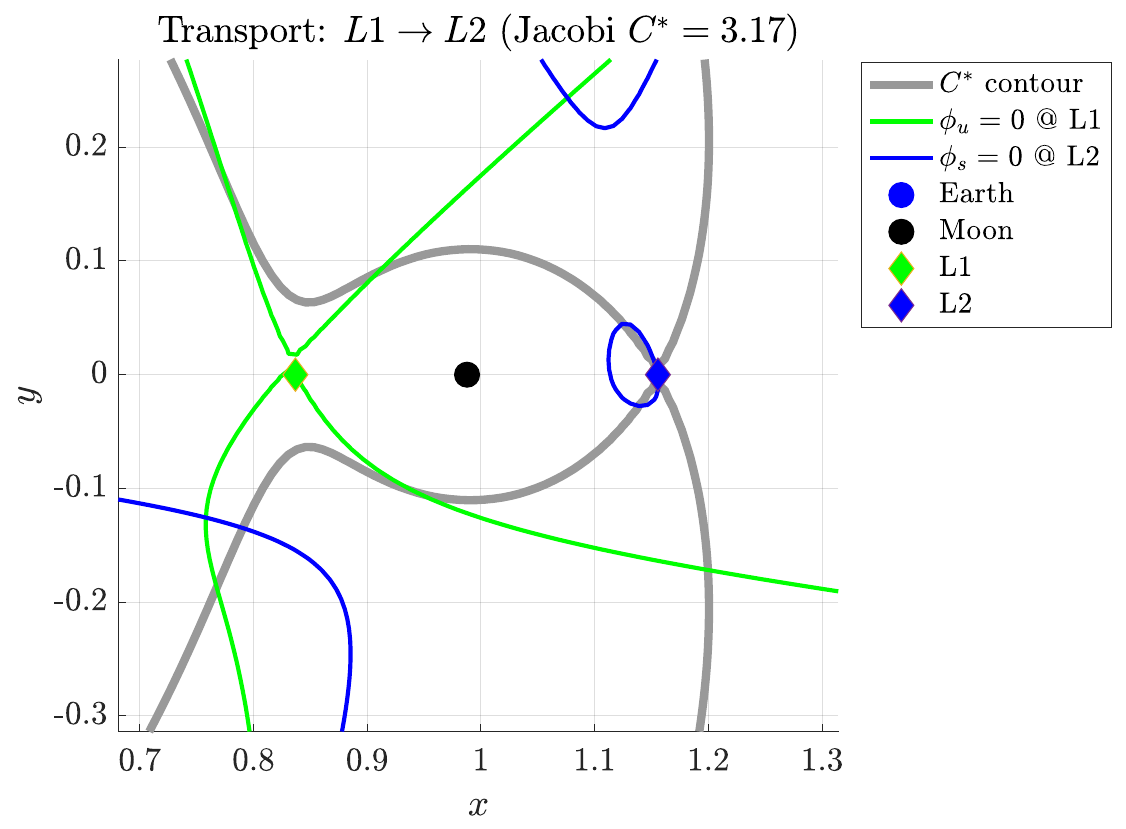}}
\caption{(a) The global alignment of zero level curves of Koopman eigenfunctions. (b) Koopman eigenfunction transport channel between L$_1$ and L$_2$.}
\label{fig_koopman_implicit}
\end{figure}

Fig. \ref{fig_koopman_implicit}(a) depicts the zero--level sets of Koopman eigenfunctions associated with the five equilibrium points of the planar CR3BP. The full phase space of the model is four--dimensional, with coordinates $(x,y,v_x,v_y)$, but here we restrict to the slice $v_x=v_y=0$, which reduces the state space to the two--dimensional configuration plane $(x,y)$. This allows the spectral geometry of the system to be visualized directly in the same coordinates where the primaries and equilibrium points reside.

\noindent\textbf{Koopman eigenfunctions.}  
For each equilibrium point $L_i$ ($i=1,\dots,5$), the Koopman operator admits an eigenfunction $\hat{\phi}_i(x,y,v_x,v_y)$ aligned with the corresponding eigendirection of the linearized dynamics. These eigenfunctions extend globally by analytic continuation along trajectories. Their zero--level sets,
\[
\Gamma_i \;=\; \big\{ (x,y) \in \mathbb{R}^2 : \hat{\phi}_i(x,y,0,0) = 0 \big\},
\]
form implicit curves in the $(x,y)$ plane. Each $\Gamma_i$ is shown in the figure as a colored contour: blue for $L_1$, orange for $L_2$, yellow for $L_3$, purple for $L_4$, and green for $L_5$. The positions of the primaries (Earth and Moon) and the Lagrange points are also marked.

\noindent\textbf{Geometric interpretation.}  
The Koopman implicit curves $\Gamma_i$ are not simply energetic contours, but rather encode the asymptotic influence of the corresponding equilibrium. On the other hand, invariant manifolds of periodic orbits and transport channels through the bottlenecks near $L_1$ and $L_2$ are naturally organized by these curves. For instance, the blue and orange contours bending into narrow channels highlight the phase space channels, through which motion can be transferred between the Earth and Moon regions. Conversely, the smooth closed shapes around $L_4$ and $L_5$ reflect their character as centers of stability.

\noindent\textbf{Spectral partitioning of phase space.}  
Whereas zero--velocity curves and Jacobi surfaces describe restrictions imposed by energy, the Koopman implicit curves partition the phase space spectrally. They identify separatrices and invariant boundaries that organize long--time dynamics. The coexistence of these operator--theoretic partitions with classical energetic structures provides a powerful framework: energy dictates where motion is possible, while Koopman eigenfunctions reveal how trajectories are asymptotically organized within those allowed regions.

\subsection{Transport Channel Between L$_1$ and L$_2$}


Fig. \ref{fig_koopman_implicit}(b) demonstrates the mechanism of transport across the Earth--Moon system.  
The unstable manifold of L$_1$ (green) threads outward through the Jacobi neck, while the stable manifold of L$_2$ (blue) enters from the opposite side.  
Although the two curves do not literally coincide, they both lie within the same Jacobi contour $C^*$, and together they bound an invariant tube in phase space.  
This tube acts as a transport channel: trajectories launched inside it will flow from the Earth region, through L$_1$, and into the Moon’s exterior region near L$_2$.  
This illustrates the essence of phase space transport: invariant manifolds, revealed by Koopman eigenfunctions, define the invisible highways that govern motion between domains of attraction. 

The results demonstrate that Koopman eigenfunctions, computed via a path--integral formulation, recover the invariant geometry of the CR3BP.  
The Jacobi contours provide global energetic constraints, while the eigenfunctions reveal local manifolds at each equilibrium.  
At collinear points (L$_1$, L$_2$, L$_3$), unstable and stable directions create transport channels, while at triangular points (L$_4$, L$_5$) center manifolds dominate, supporting oscillatory but non-transporting motion.  
The combined picture highlights how operator-theoretic tools uncover the fine structure of transport within classical celestial mechanics.

\section{Conclusions and Future Work}
In this work, we demonstrated the usefulness of the Koopman operator framework for the analysis of CR3BP dynamics. In particular, we showed that the zero-level curves of the principal Koopman eigenfunctions encode the geometry of the phase space and how one can identify the various invariant manifolds using these eigenfunctions. We further demonstrated the efficacy of the proposed approach by reconstructing the global characteristics from the eigenfunctions at the five libration points. Finally, we presented preliminary results on how this framework can be used to identify spatial highways for effective transport. The applicability of Koopman framework in the space dynamics settings opens up a whole new direction of data-driven discovery and control design, where one can leverage the linearity of the operator to analyze and design control strategies for the nonlinear CR3BP. Furthermore, again because of linearity, uncertainty analysis promises to be much simpler, and we plan to investigate these various directions in the future.

\bibliographystyle{unsrt}
\bibliography{references,ref,ref1,ref2,raktim}

\end{document}